\documentclass[]{article}

\usepackage{makeidx}         % allows index generation
\usepackage{graphicx}        % standard LaTeX graphics tool
\usepackage{multicol}        % used for the two-column index
\usepackage[bottom]{footmisc}% places footnotes at page bottom

\newcommand{\uw}{{\bf{u_w}}}
\newcommand{\uwn}{\bf{u^n_w}}
\newcommand{\uwo}{{\bf{u^{n-1}_w}}}
\newcommand{\norm}[1]{\left\lVert#1\right\rVert}
\newcommand{\ez}{{\bf{e_z}}}

\usepackage{wrapfig}

\usepackage{newtxtext}       % 
\usepackage{newtxmath}       % selects Times Roman as basic font

\usepackage{graphicx,wrapfig,lipsum}
\usepackage{amsmath}
\usepackage{bm}
\usepackage{mathtools}
\usepackage{authblk}

\makeindex             % used for the subject index
\begin{document}

\title{An efficient numerical scheme for fully coupled flow and reactive transport in variably saturated porous media including dynamic capillary effects\thanks{Project founded by VISTA, a collaboration between the Norwegian Academy of Science and Letters and Equinor.}}

\author[1]{Davide Illiano}
\author[1,2]{Iuliu Sorin Pop}
\author[1]{Florin Adrian Radu}

\affil[1]{University of Bergen}
\affil[2]{Hasselt University}
\date{}
\maketitle

\abstract{In this paper, we study a model for the transport of an external component, e.g., a surfactant, in variably saturated porous media. We discretize the model in time and space by combining a backward Euler method with the linear Galerkin finite elements. The Newton method and the L-Scheme are employed for the linearization and the performance of these schemes is studied numerically. A special focus is set on the effects of dynamic capillarity on the transport equation.}

\section{Introduction}

In this work, we concentrate on efficiently solving reactive transport models in saturated/unsaturated porous media \cite{Knabner1987, Pop2004}. Such media are observable in the section of the soil closer to the surface where, in the upper part of the domain, we have a coexistence of both water and air phases while, below the water table, the soil becomes fully saturated.

In particular, our model includes dynamic capillarity effects. The capillary pressure is commonly defined as the difference between the pressures of the two phases, in our case, the air and the water. Note that, in the Richards model, the air pressure is set to be equal to zero.

Typically, the capillary pressure is assumed to be a nonlinear decreasing function depending on the water saturation. However, numerous studies are showing that such formulation is often too simplistic and that dynamic effects, due to the changes in time of the water phase, should also be included \cite{DiCarlo2004, Fucik2010, Hassanizadeh2002, Shubao2012, Zhuang2019}.
Based on this, we consider here the system:
\begin{equation}\label{nonstandard}
\begin{split}
\partial_t \theta - \nabla \cdot \big(K(\theta,\Psi)(\nabla\Psi+\ez)\big)\ &=\ \mathbb{S}_1, \\
\Psi + p_{cap}(\theta, c)\ &=\ \tau(\theta)\partial_t \theta,\\
\partial_t (\theta c) - \nabla \cdot (D\nabla c - \uw c) + R(c)\ &=\ \mathbb{S}_2.\\
\end{split}
\end{equation}
The first equation is the Richards equation, whereas the second is an ordinary differential equation used to include the non-equilibrium effects in the capillary pressure/water content relation. Equilibrium models are obtained for $\tau = 0$. Furthermore, the third equation is the reactive transport equation. Here, $\theta$ is the water content, $\Psi$ the pressure head, $c$ the concentration of the chemical component, $K$ the conductivity, $\ez$ the unit vector in the direction opposite to gravity, $D$ the diffusion/dispersion coefficient, $\uw$ the water flux, $R(c)$ the reaction term and finally $\mathbb{S}_1$ and $\mathbb{S}_2$ are any source terms or external forces involved in the process. Note that $\uw := - K(\theta,\Psi)(\nabla\Psi+\ez)$ where $K$ is a nonlinear function depending on $\theta$ and $\Psi$. In the van Genuchten model \cite{vanGenuchten1980} one has:
\begin{equation} \label{K}
K(\theta,\Psi) = \begin{cases} K_s \theta^{\frac{1}{2}} \left[1-\left(1-\theta^{\frac{n}{n-1}}\right)^{\frac{n-1}{n}} \right]^2, &\Psi \leq 0 \\
K_s, &\Psi>0.
\end{cases}
\end{equation}
$K_s$ is the saturated conductivity and $n$ is a soil dependent parameter.

The system \eqref{nonstandard} is completed by boundary conditions for $\Psi$ and $c$, and initial conditions for $\theta$ and $c$.

The rest of the paper is organized as follows: in Section \ref{Section2} the equations are discretized and linearized. Section \ref{Section3} includes a numerical example, based on the literature \cite{Haverkamp1977}, which allows us to compare the different numerical schemes. Finally, Section \ref{Section4} will conclude this paper with our final remarks. 
 
\section{The Numerical Schemes}\label{Section2}

Applying a Euler implicit time-stepping to \eqref{nonstandard} gives a sequence of time discrete nonlinear equations. To solve them we apply different linearization schemes: the Newton method, the L-Scheme and a combination of the two \cite{Illiano2019,List2016}. They are here compared thanks to a numerical example inspired by reactive models.
 
The equations in (\ref{nonstandard}) are fully coupled due to the double dependency of the capillary pressure of both the water content $\theta$ and the concentration $c$. In general, $p_{cap}$ is a function of only $\theta$, e.g., $p_{cap}:=1/\alpha(\theta^{-1/m}-1)^{1/n}$ as presented in \cite{vanGenuchten1980}. Anyhow, it has been observed \cite{Smith1999} that, if an external component is involved, the surface tension becomes a function of the concentration $c$ and thus, the capillary pressure itself is influenced by this, i.e. $p_{cap} := p_{cap}(\theta,c)$.

In the following, we use the standard notations of functional analysis. The domain $\Omega \subset \mathbb{R}^d$, $d = 1,2$ or $3$, is bounded, open and has a Lipschitz continuous boundary $\partial \Omega$. We denote by $L^2(\Omega)$ the space of real-valued, square-integrable functions defined on $\Omega$ and  $H^1(\Omega)$ its subspace containing the functions having also the first order weak partial derivatives in $L^2(\Omega)$. $H_0^1(\Omega)$ is the space of functions belonging to $H^1(\Omega)$, having zero values on the boundary $\partial \Omega$. We denote by $< \cdot, \cdot > $ the $L^2(\Omega)$ scalar product and by $\norm{\cdot}$ the associated norm. Finally, assume that $K$ is continuous and increasing, $p_{cap}\in C^1\big((0,1], [0,\infty)\big)$ is  decreasing and $\tau \in C^1\big((0,1], [0,\infty)\big)$ .

We now combine the backward Euler method with linear Galerkin finite elements for the discretization of the problem (\ref{nonstandard}). Let $N \in \mathbb{N}$ be a strictly positive natural number, define the time step size $\Delta t \ = \ T/N$ and $t_n\ =\ n\Delta t\ (n\in{1, 2, \dots, N})$. Furthermore, $\mathbb{T}_h$ is a regular decomposition of $\Omega$, $\overline{\Omega} = \underset{\mathbb{T} \in \mathbb{T}_h}{\cup} \mathbb{T}$, into $d$-dimensional simplices, with $h$ denoting the maximal mesh diameter. The finite element space $V_h \subset H_0^1(\Omega)$ is defined by
\begin{equation}
V_h := \{ v_h\in H_0^1(\Omega)\ s.t.\ v_{h|\mathbb{T}} \in \mathbb{P}_1(\mathbb{T}), \ \mathbb{T}\in \mathbb{T}_h\},
\end{equation}
where $\mathbb{P}_1(\mathbb{T})$ denotes the space of the afine polynomials on $\mathbb{T}$.

The fully discrete Galerkin formulation of the system (\ref{nonstandard}) can be written as:

\textbf{Problem P(n)} Let $n \ge 1$ be fixed. Given $\Psi^{n-1}_h, \theta^{n-1}_h, c^{n-1}_h \in V_h$, find $\Psi^{n}_h, \theta^{n}_h, c^{n}_h \in V_h$ such that there holds
\begin{equation}\label{richards1}
<\theta^{n}_h - \theta^{n-1}_h, v_{1,h} > + \Delta t <K(\theta^{n}_h,\Psi^{n}_h)(\nabla\Psi^{n}_h+ {\ez}), \nabla v_{1,h} >\ =\ \Delta t <\mathbb{S}_1, v_{1,h} > ,
\end{equation}
\begin{equation}\label{psi1}
\Delta t<\Psi^{n}_h , v_{2,h}>  + \Delta t <p_{cap}(\theta^{n}_h,c^n_h), v_{2,h}>\ =\ <\tau(\theta^{n}_h) (\theta^{n}_h - \theta^{n-1}_h), v_{2,h}>, 
\end{equation}
and 
\begin{equation}\label{tra1}
\begin{split}
< \theta^{n}_h c^{n}_h  - \theta^{n-1}_h c^{n-1}_h, v_{3,h} >+ \Delta t  <D \nabla c^{n}_h + \uwo c^{n}_h, \nabla v_{3,h}>\\
  +  \Delta t  <R(c^{n}_h),v_{3,h} >\ =\ \Delta t <\mathbb{S}_2, v_{3,h}>,
\end{split}
\end{equation}
for all $v_{1,h},v_{2,h},v_{3,h} \in V_h$.

\emph{\textbf{Remark 1}} We use $\uwo :=- K(\theta^{n-1}_h,\Psi^{n-1}_h)(\nabla\Psi^{n-1}_h+\ez)$ for the convective term in the transport equation, for simplicity reasons. Nevertheless, all the simulations presented in this paper have also been performed with ${\uwn} := - K(\theta^{n}_h,\Psi^{n}_h)(\nabla\Psi^{n}_h+\ez)$ instead of $\uwo$ and the results were almost identical.

In the following, we propose different solving strategies for the system of equations presented above. These strategies are built on the ones discussed in \cite{Illiano2019}, extending them to the case of dynamic capillary pressure ($\tau(\theta) \neq 0$). They are either a monolithic solver of the full system, or a splitting approach obtained by solving first the flow component, using a previously computed concentration, then updating the transport equation, using the newly computed pressure and water content. In both cases, one has to iterate. Each iteration requires solving a non-linear problem, for which, either the Newton methods or the L-Scheme \cite{Illiano2019, List2016, Pop2004} are considered. These strategies are then named: monolithic-Newton scheme (MON-Newton), monolithic-L-Scheme (MON-LS), nonlinear splitting-Newton (NonLinS-Newton) and nonlinear splitting-L-Scheme (NonLinS-LS).

The index $j$ denotes the iteration index. As a rule, the iterations start with the solution obtained at the previous time step, for example $\Psi^{n,1} := \Psi^{n-1}$. This is not necessary for the L-Scheme, which is globally convergent, but it appears to be a natural choice.

\subsection{The monolithic Newton method (MON-NEWTON)}

The Newton method is a well-known linearization scheme, which is quadratic but only locally convergent. Applying the monolithic Newton method to \eqref{richards1}-\eqref{tra1} leads to

\textbf{Problem MN(n,j+1)} Let $\Psi^{n-1}_h, \theta^{n-1}_h, c^{n-1},\Psi^{n,j}_h, \theta^{n,j}_h c^{n,j}_h\in V_h$ be given, find $\Psi^{n,j+1}_h, \theta^{n,j+1}_h, c^{n,j+1}_h \in V_h$ such that
\begin{equation}\label{Newtonrichards}
\begin{split}
<\theta^{n,j+1}_h - \theta^{n-1}_h, v_{1,h} > +\Delta t <K(\theta^{n,j}_h,\Psi^{n,j}_h) (\nabla(\Psi^{n,j+1}_h)
+ \ez), \nabla v_{1,h} >  \\
+\Delta t <\partial_\theta K(\theta^{n,j}_h,\Psi^{n,j}_h) (\nabla(\Psi^{n,j}_h)+\ez)(\theta^{n,j+1}_h-\theta^{n,j}_h), \nabla v_{1,h} >\\
+\Delta t <\partial_\Psi K(\theta^{n,j}_h,\Psi^{n,j}_h) (\nabla(\Psi^{n,j}_h)+\ez)(\Psi^{n,j+1}_h-\Psi^{n,j}_h), \nabla v_{1,h} >\\
=\ \Delta t <\mathbb{S}_1,v_{1,h} >,
\end{split}
\end{equation}
\begin{equation}\label{psiNewton}
\begin{split}
\Delta t<\Psi^{n,j+1}_h , v_{2,h}> +\Delta t <p_{cap}(\theta^{n,j}_h,c^{n,j}_h), v_{2,h}>  \\
+ \Delta t <\partial_\theta p_{cap}(\theta^{n,j}_h,c^{n,j}_h) (\theta^{n,j+1}_h - \theta^{n,j}_h), v_{2,h}> + \Delta t <\partial_c p_{cap}(\theta^{n,j}_h,c^{n,j}_h)\\
 (c^{n,j+1}_h - c^{n,j}_h), v_{2,h}>
 =\ <\tau(\theta^{n,j}_h) (\theta^{n,j+1}_h - \theta^{n-1}_h), v_{2,h}>\\
 +<\partial_\theta\tau(\theta^{n,j}_h) (\theta^{n,j}_h - \theta^{n-1}_h) (\theta^{n,j+1}_h - \theta^{n,j}_h), v_{2,h}>,\\
\end{split}
\end{equation}
and 
\begin{equation}\label{Newtontransp}
\begin{split}
< \theta^{n,j}_h c^{n,j+1}_h - \theta^{n-1}_h c^{n-1}_h, v_{3,h} > 
  + \Delta t < D\nabla c^{n,j+1}_h + 
  \uwo c^{n,j+1}_h, \nabla v_{3,h}> \\
  + \Delta t<R(c^{n,j}_h), v_{3,h}>  + \Delta t < \partial_c R(c^{n,j}_h)  (c^{n,j+1}_h - c^{n,j}_h)> \\
   =\ \Delta t <\mathbb{S}_2,v_{3,h}>,
\end{split}
\end{equation}
hold true for all $v_{1,h}, v_{2,h}, v_{3,h} \in V_h$.

\subsection{The monolithic $L$-scheme (MON-LS)}

The monolithic $L$-scheme for solving \eqref{richards1}-\eqref{tra1} reads as

\textbf{Problem ML(n,j+1)} Let $\Psi^{n-1}_h, \theta^{n-1}_h, c^{n-1},\Psi^{n,j}_h, \theta^{n,j}_h c^{n,j}_h\in V_h$ be given, \\$L_1^\Psi,L_1^\theta,L_2,L_3 > 0$, big enough. 
\\Find $\Psi^{n,j+1}_h, \theta^{n,j+1}_h, c^{n,j+1}_h \in V_h$ such that
\begin{equation}\label{Lrichards}
\begin{split}
<\theta^{n,j+1}_h - \theta^{n-1}_h, v_{1,h} >  +
\Delta t <K(\theta^{n,j}_h,\Psi^{n,j}_h) (\nabla(\Psi^{n,j+1}_h)+\ez),\nabla v_{1,h} > \\
+ \Delta t <L_1^\Psi (\Psi^{n,j+1}_h -\Psi^{n,j}_h) ,\nabla v_{1,h} >  + \Delta t <L_1^\theta (\theta^{n,j+1}_h -\theta^{n,j}_h) ,\nabla v_{1,h} >\\
  =\ \Delta t <\mathbb{S}_1,v_{1,h} >,
\end{split}
\end{equation}
\begin{equation}\label{Lpsi}
\begin{split}
\Delta t<\Psi^{n,j+1}_h , v_{2,h}>\ =\ -\Delta t <p_{cap}(\theta^{n,j}_h,c^{n,j}_h), v_{2,h}>\\
 + <\tau(\theta^{n,j}_h) (\theta^{n,j+1}_h - \theta^{n-1}_h), v_{2,h}> 
+ < L_2 (\theta^{n,j+1}_h - \theta^{n,j}_h), v_{2,h}> \\
\end{split}
\end{equation}
and 
\begin{equation}\label{Ltransp}
\begin{split}
< \theta^{n,j}_hc^{n,j+1}_h -  \theta^{n-1}_h c^{n-1}_h, v_{3,h} > 
 +\Delta t< D\nabla c^{n,j+1}_h + \uwo c^{n,j+1}_h, \nabla v_{3,h}>  \\
 +\Delta t <R(c^{n,j}_h), v_{3,h}>  +  <L_3(c^{n,j+1}_h-c^{n,j}_h),  v_{3,h}>
\ =\ \Delta t<\mathbb{S}_3,v_{3,h}>,
\end{split}
\end{equation}
hold true for all $v_{1,h}, v_{2,h}, v_{3,h} \in V_h$.  

The L-Scheme does not involve the computations of derivatives, the linear systems to be solved within each iteration are better conditioned, compared to the ones given by the Newton method \cite{Illiano2019, List2016}, and it is globally (linearly) convergent.
The convergence of the scheme has been proved, for the equilibrium model ($\tau(\theta)=0$) in \cite{Illiano2019}, and can be easily extended to the non-equilibrium formulation given by the system (\ref{Lrichards})-(\ref{Ltransp}).

\subsection{The splitting approach (NonLinS)}

The splitting approach for solving \eqref{richards1}-\eqref{tra1} reads as

\textbf{Problem S(n,j+1)} Let $\Psi^{n-1}_h, \theta^{n-1}, c^{n-1},\Psi^{n,j}_h, \theta^{n,j}_h, c^{n,j}_h\in V_h$ be given, find $\Psi^{n,j+1}_h, \theta^{n,j+1}_h \in V_h$ such that
\begin{equation}\label{NonLinrichards}
\begin{split}
<\theta^{n,j+1}_h - \theta^{n-1}_h, v_{1,h} > +\Delta t <K(\theta^{n,j+1}_h,\Psi^{n,j+1}_h) (\nabla(\Psi^{n,j+1}_h)+ \ez), \nabla v_{1,h} >\\
=\ \Delta t <\mathbb{S}_1,v_{1,h} >,
\end{split}
\end{equation}
\begin{equation}\label{NonLinpsi}
\begin{split}
\Delta t<\Psi^{n,j+1}_h , v_{2,h}> + \Delta t <p_{cap}(\theta^{n,j+1}_h,c^{n,j}_h), v_{2,h}>\ \\
=\ <\tau(\theta^{n,j+1}_h) (\theta^{n,j+1}_h - \theta^{n-1}_h), v_{2,h}> ,
\end{split}
\end{equation}
hold true for all $v_{1,h},v_{2,h} \in V_h$.

Then, with  $\Psi^{n,j+1}_h$ and $\theta^{n,j+1}_h$ obtained from the equations above, find $c^{n,j+1}_h \in V_h$ such that
\begin{equation}\label{NonLintransp}
\begin{split}
< \theta^{n,j+1}_h c^{n,j+1}_h - \theta^{n-1}_hc^{n-1}_h, v_{3,h} >  
 + \Delta t < D\nabla c^{n,j+1}_h + \uwo c^{n,j+1}_h, \nabla v_{3,h}> \\
 +\Delta t <R(c^{n,j+1}_h), v_{3,h}>\ =\ \Delta t <\mathbb{S}_2,v_{3,h}>,
\end{split}
\end{equation}
holds true for all $v_{3,h} \in V_h$.

The three equations above can be then linearised using either the Newton method (NonLinS-Newton) or the L-Scheme (NonLinS-LScheme).

\subsection{The mixed linearization scheme}
It has been already observed, for a different set of equations \cite{List2016}, that combining the Newton method and the L-Scheme can improve the convergence of the scheme.
The Newton method is quadratically but only locally convergent and it can produce badly conditioned linearized systems. Moreover, the time step is subject to severe restrictions for guaranteeing the convergence of the scheme, and this has also been observed in numerical examples \cite{Cao2015, Illiano2019, List2016}.

Contrarily, the L-Scheme is globally convergent and the linear systems to be solved within each iteration are better conditioned, however, it has only a linear rate of convergence. 

The mixed formulation, obtained combining the two schemes, appears to be the best approach and shows practically both global and quadratic convergence. 
The Newton method commonly fails to converge, if the initial guess is too far from the actual solution. Since this guess is usually the solution at the previous time, this can force restriction on the time step. Instead of reducing the time step one can obtain a better approximation of the initial guess, for the Newton method, by performing few L-Scheme iterations. In the numerical simulation here presented, up to 5 iterations were sufficient to reach a good initial guess for the Newton iteration, which ensured its convergence.

\section{Numerical examples}\label{Section3}
In this section, we use a benchmark problem, from \cite{Haverkamp1977}, to compare the different linearization schemes and solving algorithms defined above. It describes the recharge of a two-dimensional underground reservoir $\Omega \subset \mathbb{R}^2$, in the interval of time $t\in (0, 3]$. The boundary of the domain and the Dirichlet boundary conditions are defined below.

\noindent\begin{minipage}{.5\linewidth}
\begin{equation*}
  \begin{split}
&\Omega = (0,2) \times (0,3), \\
&\Gamma_{D_1} = \{(x,y) \in \partial \Omega | x \in [0,1] \wedge y = 3\}, \\ 
&\Gamma_{D_2} = \{(x,y) \in \partial \Omega | x = 2 \wedge y \in [0,1]\}, \\
&\Gamma_D = \Gamma_{D_1} \cup \Gamma_{D_2}, \\
&\Gamma_N = \partial \Omega \setminus \Gamma_D,
\end{split}
\end{equation*}
\end{minipage}%
\begin{minipage}{.5\linewidth}
\begin{equation*}
  \begin{split}
&\Psi(x,y,t) = \begin{cases} 
-2 + 2.2 *t,\  \text{on}\ \Gamma_{D_1}, t\leq 1 \\
0.2,\ \text{on}\ \Gamma_{D_1},  t > 1\\
1-y,\ \text{on}\ \Gamma_{D_2},
\end{cases}\\
&c(x,y,t) = \begin{cases} 
1,\  \text{on}\ \Gamma_{D_1}, t\leq 1 \\
0,\ \text{on}\ \Gamma_{D_1},  t > 1\\
3-y,\ \text{on}\ \Gamma_{D_2} \cup \Gamma_N.
\end{cases}\\
\end{split}
\end{equation*}
\end{minipage}

Furthermore, no flow conditions are imposed on $\Gamma_N$. The initial conditions are given by $\Psi(x,y,0):=1-y, c(x,y,0):=3-y$ and $\theta(x,y,0):=0.39$. The capillary pressure is defined as $p_{cap}(\theta,c):= (1 - \theta)^{2.5} + 0.1*c$, the conductivity is given by (\ref{K}) and $\tau(\theta) = 1$. Finally, the parameters implemented are: $K_s = 1$, $L_1^\Psi,L_1^\theta,L_2 = 0.01$, $L_3 = 0.1$ and the iterations stop whenever all the error norms, $\norm{\Psi^{n,j+1} - \Psi^{n,j}},\norm{\theta^{n,j+1} - \theta^{n,j}}$ and $\norm{c^{n,j+1} - c^{n,j}}$, are below $10^{-6} $.

\begin{figure}[h]
\centering
\includegraphics[scale=0.5]{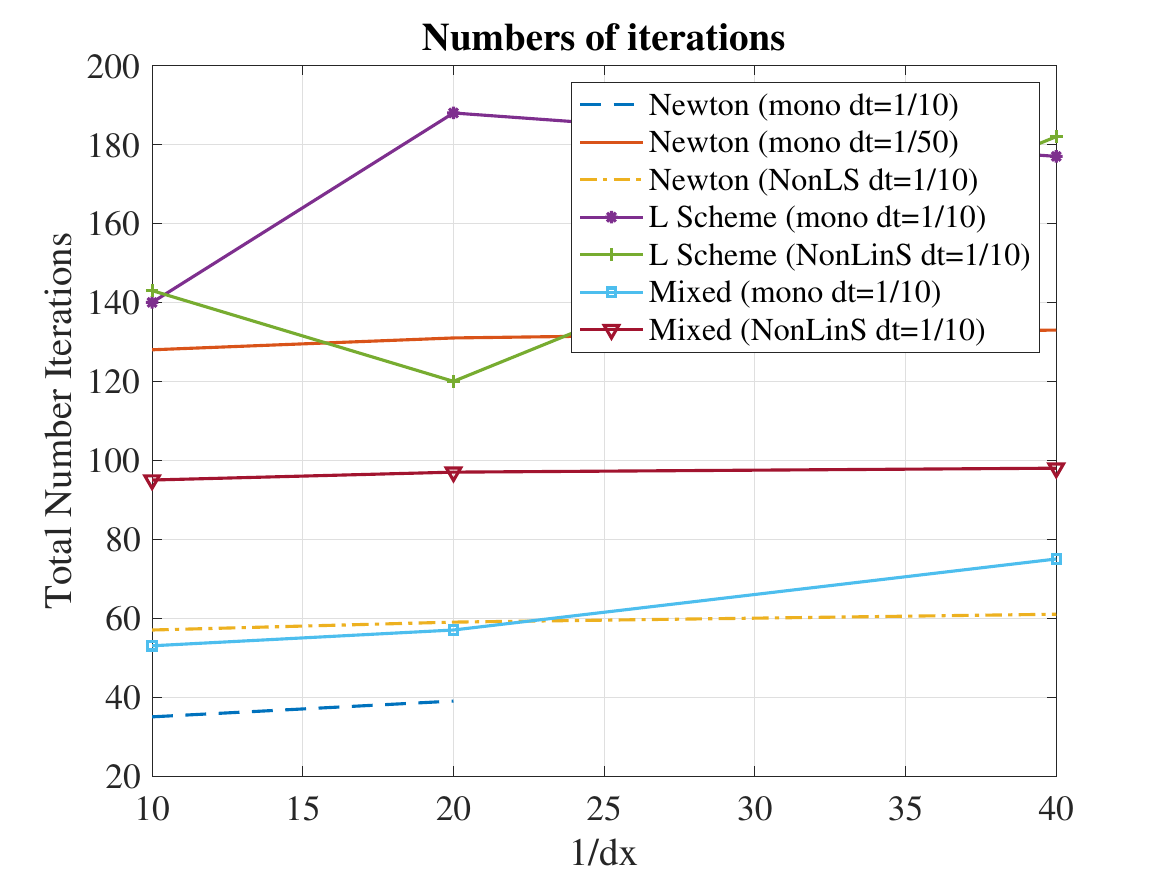}
\caption{Total numbers of iterations for different solvers}\label{Figure1it}
\end{figure}

We performed the simulations using regular meshes, consisting of squares, with sides $dx = \{1/10,1/20,1/40\}$. We considered two fixed time steps $\Delta t = 1/10$ and $\Delta t = 1/50$. 

In Figure \ref{Figure1it}, we can observe the total numbers of iterations required by the different linearization schemes and solving algorithms. Next to the name of each scheme we report, between parenthesis, which time step $\Delta t$ has been used.

We can observe, as the Newton method in the monolithic formulation, converges only for coarse meshes, for $\Delta t=1/10$. For the smaller time step, $\Delta t = 1/50$, it converges for all of the tested meshes.

The L Scheme converges for both time steps, but, since it is linearly convergent, for $\Delta t = 1/50$ would require more iterations than the Newton method. 

The results obtained thanks to the mixed formulation are particularly interesting. We can observe that this scheme, both in the monolithic and splitting formulation, converges for all the tested meshes also in case of a large time step. Moreover, thanks to the Newton iterations, it appears to be faster than the classical L Scheme. It is as robust as the L Scheme and as fast as the Newton method. For more details regarding the mixed scheme, we refer to \cite{List2016}.

\section{Conclusions}\label{Section4}

In this paper, we considered multiphase flow coupled with a one-component reactive transport in variably saturated porous media, including also the dynamic effects in the capillary pressure. The resulting model is nonlinear and for this reason, three different linearization schemes are investigated: the L-Scheme, the Newton method and a combination of the two. We also studied both monolithic solvers and splitting ones. 

The tests show that, for this particular set of equations, the best linearization scheme is the one obtained combining the Newton method and the L-Scheme. Such scheme appears to be both quadratically and globally convergent. 

\vspace{2cm}
\emph{Acknowledgments}
\\The research of D. Illiano was funded by VISTA, a collaboration between the Norwegian Academy of Science and Letters and Equinor, project number 6367, project name: adaptive model and solver simulation of enhanced oil recovery. The research
of I.S. Pop was supported by the Research Foundation-Flanders (FWO), Belgium through the Odysseus programme (project G0G1316N) and Equinor through the Akademia grant. 
\\We thank the members of the \emph{Sintef} research group and in particular to Dr. Olav Moyner for the assistance with the implementation of the numerical examples in \emph{MRST}, the toolbox based on Matlab developed at \emph{Sintef} itself.


\begin{thebibliography}{99.}%
% Use the following syntax and markup for your references if 
% the subject of your book is from the field 
% "Mathematics, Physics, Statistics, Computer Science"
%
%

\bibitem{Cao2015}
Cao, X., Pop, I.S.: Uniqueness of weak solutions for a pseudo-parabolic equation modeling two phase flow in porous media,
Applied Mathematics Letters, Volume 46, Pages 25-30, (2015).

\bibitem{DiCarlo2004}
Di Carlo, D.: Experimental measurements of saturation overshoot on infiltration, Water Resources Research, Volume 40, Issue 4, (2004).

\bibitem{Fucik2010}
Fucik, R., Mikyska, J., Sakaki, T., Benes, M., Illangasekare, T.H.: Significance of Dynamic Effect in Capillarity during Drainage Experiments in Layered Porous Media, Vadose Zone Journal 9, Volume 3, (2010).

\bibitem{vanGenuchten1980}
van Genuchten, M.: A Closed-form Equation for Predicting the Hydraulic Conductivity of Unsaturated Soils, Soil Science Society of America Journal, Volume 44, Issue 5, Pages 892-898, (1980).

\bibitem{Hassanizadeh2002}
Hassanizadeh, S.M., Celia, M.A., Dahle, H.K.: Dynamic Effect in the Capillary Pressure Saturation Relationship and its Impacts on Unsaturated Flow, Vadose Zone Journal, Volume 1, Issue 1, Pages 38-57, (2002).
 
 \bibitem{Haverkamp1977}
Haverkamp, R., Vauclin, M., Touma, J., Wierenga, P.J., Vachaud, G.: A comparison of numerical simulation models for one-dimensional infiltration, Soil Science Society of America Journal, Volume 41, Pages 285-294, (1977).

\bibitem{Illiano2019} Illiano, D., Pop, I.S., Radu, F.A.: Iterative schemes for surfactant transport in porous media, arXiv preprint arXiv:1906.00224, (2019).

\bibitem{Knabner1987}
Knabner, P.: Finite element simulation of saturated-unsaturated flow through porous media, LSSC 7, Pages 83-93, (1987).

\bibitem{List2016} List, F., Radu, F.A.: A study on iterative methods for solving Richards' equation, Computational Geoscience, Volume 20, Issue 2, Pages 341-353, (2016).

\bibitem{Pop2004}
Pop, I.S., Radu, F.A., Knabner, P.: Mixed finite elements for the  Richards' equation: linearization procedure, Journal of computational and applied mathematics, Volume 168, Issue 1, Pages 365-373, (2004).
 
\bibitem{Shubao2012}
Shubao, T., Lei, G., Shun-li, H., Yang, L.: Dynamic effect of capillary pressure in low permeability reservoirs, Petroleum Exploration and Development, Volume 39, Issue 3, Pages 405-411, (2012).

\bibitem{Smith1999}
Smith, J., Gillham, R.: Effects of solute concentration-dependent surface tension on unsaturated flow: Laboratory sand column experiments, Water Resource Research, Volume 35, Issue 4, Pages 973-982, (1999).

\bibitem{Zhuang2019}
Zhuang, L., van Duijn, C.J.,  Hassanizadeh, S.M.: The effect of dynamic capillarity in modeling saturation overshoot during infiltration, Vadose Zone Journal, Volume 18, Issue 1, Pages 1-14, (2019).

\end{thebibliography}
\end{document}